\documentclass[12pt]{article}
\usepackage{amsfonts}
\usepackage{amsmath}
\usepackage{amssymb,amsmath}
\def\fc{{\textbf{\textit c}}}
\def\D{\Delta}
\def\a{\alpha}

\def\b{\beta}
\def\SM{\!\setminus\!}

\def\UU{{\mathcal U}}
\def\LL{{\mathcal L}}

\def\Der{{\rm Der}}
\def\Inn{{\rm Inn}}
\def\Ker{{\rm Ker}}
\def\Im{{\rm Im}}

\def\cl{\centerline}
\def\rar{\rightarrow}

\def\vs{\vspace*}

\def\ni{\noindent}
\def\VV{\mathcal {V}}
\def\Z{\mathbb{Z}}

\def\F{\mathbb{F}}
\def\QED{\hfill$\Box$}

\numberwithin{equation}{section}
\newtheorem{theo}{Theorem}[section]
\newtheorem{defi}[theo]{Definition}

\newtheorem{lemm}[theo]{Lemma}
\newtheorem{prop}[theo]{Proposition}

\begin{document}

\cl{{\large\bf Lie bialgebra structures on extended}} \cl{{\large\bf
Schr\"{o}dinger-Virasoro Lie algebra}\footnote {Supported by
NSF grants 10825101 of China China \\[2pt]\indent $^{*}$Corresponding author:
lmyuan@mail.ustc.edu.cn}}\vs{6pt}

\cl{ Lamei Yuan$^{\,*}$, Yongping Wu$^{\,\dag}$, Ying Xu$^{\,*}$}
 \cl{\small $^{*}$Department of Mathematics,
University of Science \!and \!Technology \!of \!China, Hefei 230026,
China} \cl{$^{\dag}$School of Mathematics and Computer Science,
Longyan College, Fujian 364111, China}
 \cl{\small E-mail: lmyuan@mail.ustc.edu.cn, wyp\_504@sohu.com, xying@mail.ustc.edu.cn}\vs{6pt}

{\small\parskip .005 truein \baselineskip 3pt \lineskip 3pt

\noindent{{\bf Abstract.} In this paper, Lie bialgebra structures on
the extended Schr\"{o}dinger-Virasoro Lie algebra $\LL$ are
classified. It is obtained  that all the Lie bialgebra structures on
$\LL$ are triangular coboundary. As a by-product, it is derived that
the first cohomology group $H^1(\LL,\LL\otimes\LL)$ is trivial.
\vs{5pt}

}

\vs{5pt}


\parskip .001 truein\baselineskip 6pt \lineskip 6pt

\vs{18pt}

\cl{\bf\S1. \
Introduction}\setcounter{section}{1}\setcounter{equation}{0}

\vs{8pt}

The Schr\"{o}dinger-Virasoro Lie algebras \cite{H1} were introduced
in the context of non-equilibrium statistical physics during the
process of investigating the free Schr\"{o}dinger equations. They
are closely related to the Schr\"{o}dinger algebra and the Virasoro
algebra, both of which play important roles in many areas of
mathematics and physics (e.g., statistical physics, integrable
system) and have been investigated in a series of papers
[\ref{HU}--\ref{LS2}, \ref{LSZ}, \ref{RU}, \ref{TZ}, \ref{ZT}]. In
order to investigate vertex representations of the
Schr\"{o}dinger-Virasoro Lie algebra, J. Unterberger introducted
(see \rm {Definition 1.5} in [\ref{U}]) a class of
infinite-dimensional Lie algebras called the extended
Schr\"{o}dinger-Virasoro Lie algebra $\LL$, which can be viewed as
an extension of the Schr\"{o}dinger-Virasoro Lie algebra by a
conformal current of weight $1$ and generated by $\{L_n,M_n, N_n,
Y_p\,|\,n \in \Z, \ p\in \Z+1/2\}$ with the following Lie brackets:
\begin{eqnarray}\label{LB}\begin{array}{lll}
&&[L_m,L_{n}]=(n-m)L_{n+m},\ \ \ \ \ [L_m,N_n]=nN_{m+n},\ \ \ \ \,
[L_m,M_n\,]=nM_{n+m},\
\\[4pt]&&[L_n,\,Y_p\,]=(p-n/2)Y_{p+n},\ \ \ \ \, [N_m,Y_p\,]=Y_{m+p},\ \ \ \ \ \ \ \,[N_m,M_n]=2M_{m+n},\ \ \\ [4pt]&&
[M_n,Y_p\,]=[N_m,N_n]=0,\ \ \ \ \ \ [M_m,M_n]=0, \ \ \ \ \ \ \ \ \ \
\ \,[\,Y_p,Y_{q}\,\,]=(q-p)M_{p+q}.\end{array}
\end{eqnarray}
Note that $\LL$ is centerless and finitely generated with a
generating set $\{L_{-2},L_{-1},L_1,L_2,N_1,Y_{1/2}\}$. Moreover, it
is $\frac12{\mathbb{Z}}$ -graded by
$$\LL=\mbox{$\bigoplus\limits_{n\in\Z}$}\LL_{n/2}=\big(\mbox{$\bigoplus\limits_{n\in Z}$}\LL_n\big)\mbox{$\bigoplus$}\big(\mbox{$\bigoplus\limits_{n\in\Z}$}\LL_{n+1/2}\big),$$
where $\LL_n=span\{L_n,M_n,N_n\}$ and
$\LL_{n+1/2}=span\{Y_{n+1/2}\}$, for all $n\in\Z$. The derivations,
central extensions and automorphisms of $\LL$ have been studied in
[\ref{GJP}].

 To search for the solutions of the Yang-Baxter quantum
equation, Drinfel'd \cite{D1} introduced the notion of Lie
bialgebras in 1983. Since then, a number of people have studied
further Lie bialgebra structures (e.g., [\ref{EK}, \ref{LSX1},
\ref{LSX2}, \ref{M1}--\ref{M3}). Witt type Lie bialgebras introduced
in \cite{T} were classified in \cite{NT}. This work has been
generalized  in \cite{SS,WSS1}. Lie bialgebra structures on
generalized Virasoro-like and Block  Lie algebras were investigated
in \cite{ LSX1, WSS2}.
 Drinfel'd \cite{D2} posed
the problem that whether or not there exists a general way to
quantilize all Lie bialgebras. Etingof and Kazhdan \cite{EK} gave a
positive answer to this problem, but there does not exist an uniform
method to realize quantilizations of all Lie bialgebras. Actually,
investigating Lie bialgebras and quantilizations is a complicated
problem. The authors in [\ref{HS}] prove that not all Lie bialgebra
structures on the Schr\"{o}dinger-Virasoro Lie algebra are
triangular coboundary. For the extended Schr\"{o}dinger-Virasoro Lie
algebra $\LL$, this is not the case. Namely,  we obtain that all Lie
bialgebra structures on $\LL$ are triangular coboundary. In
particular, we derive that the first cohomology group
$H^1(\LL,\LL\otimes\LL)$ is trivial.

\vs{8pt}

 \cl{\bf\S2. \ Preliminaries}\setcounter{section}{2}
\setcounter{theo}{0}\setcounter{equation}{0}

\vs{8pt}

 Throughout the paper, $\F$ denotes a filed with characteristic
 zero.
 All vector spaces and tensor products are over $\F$. Let $\Z_+$ (resp. $\Z_{>0}$) be the
 set of all nonnegative (resp. positive) integers and $\Z^*$
 be the set of all nonzero elements of $\Z$.\par
Let $L$ be a vector space, $\xi$ the {\it cyclic map} of $L\otimes
L\otimes L$, namely, $ \xi (x_{1} \otimes x_{2} \otimes x_{3})
=x_{2} \otimes x_{3} \otimes x_{1}$ for $x_1,x_2,x_3\in L,$ and
$\tau$ the {\it twist map} of $L\otimes L$, i.e., $\tau(x\otimes y)=
y \otimes x$ for $x,y\in L$.
 A {\it Lie algebra} is
a pair $(L,\delta)$, where $\delta :L\otimes L\rar L$ is a bilinear
map with the conditions:
\begin{eqnarray*}
&&\Ker(\,{\rm Id}-\tau) \subset \Ker\,\delta,\ \ \ \delta \cdot
(\,{\rm Id} \otimes \delta ) \cdot (\,{\rm Id} + \xi +\xi^{2}) =0 :
L \otimes L\otimes L\rar
 L,
\end{eqnarray*}
where ${\rm Id}$ is the identity map. Dually, a {\it Lie coalgebra}
is a pair $(L,\D)$ with a linear map $\D: L\to L\otimes L$
satisfying:
\begin{eqnarray}\label{cLie-s-s}
&&\Im\,\D \subset \Im(\,{\rm Id}- \tau),\ \ \ (\,{\rm Id} + \xi
+\xi^{2}) \cdot (\,{\rm Id}\otimes \D) \cdot \D =0: L\to L\otimes
L\otimes L.
\end{eqnarray}
A {\it Lie bialgebra} is a triple $( L,\delta,\D)$ such that $( L,
\delta)$ is a Lie algebra, $( L,\D)$ is a Lie coalgebra, and the
following compatible condition holds:
\begin{eqnarray}
\label{tr}&&\D\delta (x\otimes y) = x \cdot \D y - y \cdot \D x,\ \
\forall\,\,x,y\in L.
\end{eqnarray}
where ``$\cdot$'' means the {\it diagonal adjoint action}, i.e.,
$x\cdot(\mbox{$\sum_{i}$}{a_{i}\otimes b_{i}})=
\mbox{$\sum_{i}$}({[x,a_{i}]\otimes b_{i}+a_{i}\otimes[x,b_{i}]}),$
and in general, $\delta (x\otimes y)=[x,y]$, for all $x,y,a_i,
b_i\in L$.

Denote by $\UU$ the universal enveloping algebra of $ L$ and  $1$
the identity element of $\UU$. For any $r =\sum_{i} {a_{i} \otimes
b_{i}}\in L\otimes L$, define $\fc(r)\in \UU \otimes \UU \otimes
\UU$ by
\begin{eqnarray*}
&&\fc(r) = [r^{12} , r^{13}] +[r^{12} , r^{23}] +[r^{13} , r^{23}],
\end{eqnarray*}
where $r^{12}=\sum_{i}{a_{i} \otimes b_{i} \otimes 1} , \ \ r^{13}=
\sum_{i}{a_{i} \otimes 1 \otimes b_{i}} , \ \ r^{23}=\sum_{i}{1
\otimes a_{i} \otimes b_{i}}$. Obviously,
\begin{eqnarray*}
\fc(r)=\mbox{$\sum\limits_{i,j}$}[a_i,a_j]\otimes b_i\otimes
b_j+\mbox{$\sum\limits_{i,j}$}a_i\otimes [b_i,a_j]\otimes b_j+
\mbox{$\sum\limits_{i,j}$}a_i\otimes a_j\otimes [b_i,b_j].
\end{eqnarray*}

\begin{defi}\label{def2}
\rm (1) A {\it coboundary Lie bialgebra} $( L, \delta, \D,r)$ is a
Lie bialgebra such that the cobracket $\Delta$ is an inner
derivation, i.e., there exists an element $r\in L\otimes L$ such
that
\begin{eqnarray}
\label{D-r}\D(x)=x\cdot r\mbox{\ \ for\ all\ \ }x\in L.
\end{eqnarray}
$\Delta$ is called a coboundary of $r$, denoted by $\Delta_r$.\par
 (2) A coboundary Lie bialgebra $( L,\delta,\D,r)$ is called
{\it triangular} if it satisfies the following {\it classical
Yang-Baxter Equation} (CYBE):
\begin{eqnarray}
\label{CYBE} \fc(r)=0.
\end{eqnarray}

(3) An element $r\in\Im(\,{\rm Id}-\tau)\subset L\otimes L$ is said
to satisfy the \textit{modified Yang-Baxter equation} (MYBE) if
\begin{eqnarray}
\label{MYBE}x\cdot \fc(r)=0,\ \,\forall\,\,x\in L.
\end{eqnarray}
\end{defi}\par
The following famous results are due to Drinfel'd [\ref{D2}],
Michaelis \cite{M1} and Taft \cite{T}, respectively. We combine them
into one theorem as follows:
\begin{theo}\label{theo}\vskip-3pt
 \begin{itemize}\parskip-2pt
\item[\rm(i)] For a Lie algebra $(L,[\cdot,\cdot])$ and $r\in \mathrm{Im}(\,{\rm Id}-\tau)\subset L\otimes
L$, the triple $ (L,[\cdot,\cdot],\Delta_r)$ is a Lie bialgebra if
and only if $r$ satisfies MYBE.

\item[\rm(ii)]  Let $L$ be a Lie algebra containing two linear
independent elements $a,b$ satisfying $[a,b]=kb$ for some nonzero
$k\in \F$, and set $r=a\otimes b-b\otimes a$. Then $r$ is a solution
of CYBE and equips $L$ with a structure of triangular coboundary Lie
bialgebra.

\item[\rm(iii)] Let $L$ be a Lie algebra and
$r\in\mathrm{Im}(\,{\rm Id}-\tau)\subset L\otimes L$. Then for any
$x\in L$,
\begin{equation*}
(\,{\rm Id}+\xi+\xi^2)\cdot(1+\Delta_r)\cdot\Delta_r(x)=x\cdot c(r).
\end{equation*}
\end{itemize}\end{theo}
\cl{\bf\S3. \ Structures of Lie bialgebra of the extended
Schr\"{o}dinger-Virasoro  Lie algebra }\setcounter{section}{3}
\setcounter{theo}{0}\setcounter{equation}{0}

\vs{8pt}

Regard $\VV=\LL\otimes\LL$ as a $\LL$-module under the adjoint
diagonal action. A linear map $D:\LL\to\VV$ is called a
\textit{derivation} if
\begin{eqnarray}
\label{deriv} D([x,y])=x\cdot D(y)-y\cdot D(x) \mbox{\ \ for\ all\ \
}x,\ y\ \in \LL.
\end{eqnarray}
If there exists some $v\in\VV$ such that $D(x)=x\cdot v$, then $D$
is called an \textit{inner derivation}. Denote by $v_{inn}$ the
inner derivation determined by $v$. Let $\Der(\LL,\VV)$ (resp.
$\Inn(\LL,\VV)$) be the set of all derivations (resp. inner
derivations). Then it is well known that
$H^1(\LL,\VV)\cong\Der(\LL,\VV)/\Inn(\LL,\VV), $ where
$H^1(\LL,\VV)$ is the {\it first cohomology group} of the Lie
algebra $\LL$ with coefficients in $\VV$.

 A derivation $D\in\Der(\LL,\VV)$ is {\it homogeneous of degree
$\alpha\in\frac12\Z$} if $D(\LL_p)\subset \VV_{\alpha+p}$ for all
$p\in\frac12\mathbb{Z}$. Let $\Der(\LL,\VV)_\alpha$ be the set of
all the homogeneous derivations of degree $\alpha$.
 For any $D\in\Der(\LL,\VV)$ and
$\alpha\in\frac12\Z$, define a linear map
$D_\alpha:\LL\rightarrow\VV$ as follows: for any $\mu\in\LL_q$ with
$q\in\frac12\mathbb{Z}$, write
$D(\mu)=\sum_{p\in\frac12\mathbb{Z}}\mu_p$ with $\mu_p\in\VV_p$,
then we set $D_\alpha(\mu)=\mu_{q+\alpha}$. Obviously, $D_\alpha\in
\Der(\LL,\VV)_\alpha$ and we have
\begin{eqnarray}\label{summable}
D=\mbox{$\sum_{\alpha\in\frac12\mathbb{Z}}D_\alpha$},
\end{eqnarray}
which holds in the sense that only finitely many $D_\alpha(u)\neq 0$
and $D(u)=\sum_{\alpha\in\frac12\mathbb{Z}}D_\alpha(u)$ for any
$u\in\LL$. Actually, for any $D\in{\rm Der}(\LL,\VV)$,
(\ref{summable}) is a finite sum, referring to [\ref{HS}] for
details.

\begin{lemm}\label{lmma1}\rm $H^1(\LL_0, \VV_{n/2})=0$ for all
$n\in \Z^*$.
\end{lemm}\vs{-6pt}
\noindent{\it Proof}.\ \ For any $D\in Der(\LL,\VV)$, we have
$D=\sum_{n\in\Z}D_{n/2}$. Suppose $n\neq 0$, then the restriction of
$D_{n/2}$ to $\LL_0$ induces a derivation from $\LL_0$ to the
$\LL_0$-module $\VV_{n/2}$. That is, $D_{n/2}|_{\LL_0}\in
\Der(\LL_0, \VV_{n/2})$. Conveniently, we denote $D_{n/2}|_{\LL_0}$
by $D_{n/2}$. Let $r=\frac{2}{n} D_{n/2}(L_0)\in \VV_{n/2}$. For any
$X_0\in \LL_0$, one has $\frac{n}{2}D_{n/2}(X_0)=L_0\cdot
D_{n/2}(X_0)=X_0\cdot D_{n/2}(L_0)$, since $[L_0, X_0]=0$. It
follows $D_{n/2}(X_0)=X_0\cdot r$, which implies $D_{n/2}$ is inner.
\QED
\begin{lemm}\label{lmma2}\rm $\rm{Hom}_{\LL_0}(\VV_{m/2}, \VV_{n/2})=0$ for all
$m\neq n$.
\end{lemm}\vs{-6pt}
\noindent{\it Proof}.\ \ Let $f\in\rm{ Hom}_{\LL_0}(\VV_{m/2},
\VV_{n/2})=0$ with $m\neq n$. One has $f([X_0,E_{m/2}])=[X_0,
f(E_{m/2})]$ for any $X_0\in \LL_0$ and $E_{m/2}\in \VV_{m/2}$. In
particular,  $f([L_0,E_{m/2}])=[L_0, f(E_{m/2})]$. That is,
$\frac{m}{2}f(E_{m/2})=\frac{n}{2}f(E_{m/2})$. It follows
$f(E_{m/2})=0$, since $m\neq n$. Consequently, $f=0$.\QED\par

Taking these two Lemmas above into account, we can immediately
derive the following result from Proposition 1.2 in \cite{F}.
\begin{prop}\label{p2}\rm $\Der(\LL, \VV)=\Der_0(\LL,
\VV)+\Inn(\LL,\VV)$.
\end{prop}
\begin{lemm}\label{lemm3}\rm
Let $\LL^{\otimes n}=\LL\otimes\cdots\otimes\LL$ be the tensor
product of $n$ copies of $\LL$, and regard $\LL^{\otimes n}$ as an
$\LL$-module under the adjoint diagonal action. Suppose
$r\in\LL^{\otimes n}$ satisfying $x\cdot r=0$, $\forall$ $x\in\LL$.
Then $r=0$.
\end{lemm}\vs{-6pt}
\noindent{\it Proof.}\ \  It is easy to see that $\LL^{\otimes n}$
is $\frac{1}{2}\Z$-graded by \par
$$\LL^{\otimes
n}_{p}=\mbox{$\sum\limits_{p_1+p_2+\cdots+p_n=p}\LL_{p_1}\otimes
\LL_{p_2}$}\otimes\cdots \otimes \LL_{p_n}, \ \ \forall \ \ p, \
p_i\in \frac{1}{2}\Z, \ \ i=1,2,\cdots n .
$$Write $r=\sum_{p\in\frac{1}{2}\Z}r_{p}$ as a finite sum with $r_p\in
\LL^{\otimes n}_p$. By hypothesis, $L_0\cdot r=0$, which implies
$r=r_0$. So $r=\sum\limits_{r_1+r_2+\cdots+r_n=0}c_{r_1,r_2,\cdots
r_n}E_{r_1}\otimes E_{r_2}\otimes\cdots E_{r_n}$ for some $
c_{r_1,r_2,\cdots r_n}\in \F$ and $E_{r_i}\in\LL_{r_i}$ with
$r_i\in\frac{1}{2}\Z$. Since $M_0\cdot r=0$ by the assumption, all
the coefficients of the terms containing $N_j$ for $j\in\Z$ are
zero, hence these terms in the sum vanish. Similarly, by $N_0\cdot
r=0$, one can kill the coefficients of the terms containing $M_j$
and $Y_{j+1/2}$ with $j\in\Z$. Now we can rewrite $r$ by
$$r=\mbox{$\sum\limits_{r_1+r_2+\cdots+r_n=0}c_{r_1,r_2,\cdots
r_n}$}L_{r_1}\otimes L_{r_2}\otimes\cdots L_{r_n} \ \mbox{ for
some}\ \ \  c_{r_1,r_2,\cdots r_n}\in \F.$$ But $M_1\cdot r=0$
forces all the coefficients $c_{r_1,r_2,\cdots r_n}$ are zero. This
proves the lemma. \QED
\begin{theo}\label{theo1}\vskip-3pt
$\Der(\LL, \VV)=\Inn(\LL, \VV)$.
\end{theo}
\noindent{\it Proof}.\ \ It suffices to show $\Der_0(\LL,
\VV)\subseteq \Inn(\LL,\VV)$ by virtue of Proposition \ref{p2}. For
any $0\neq D\in \Der_0(\LL,\VV)$, we shall prove that the zero
derivation is obtained after a number of steps in each of which $D$
is replaced by $D-u_{\rm inn}$ for some $u\in \VV_0$. This will be
done by a little bit complicated calculations. For clarity, we
divide them into three claims.
\par
{\bf {Claim 1.}} $D(L_0)=0$. \par In fact, for any $X_p\in \LL$ with
$p\in\frac{1}{2}\Z$, applying $D$ to $[L_0, X_p]=pX_p$, one has
$X_p\cdot D(L_0)=0$. Then it follows from Lemma \ref{lemm3} that
$D(L_0)=0$.\par

{\bf {Claim 2.}} $D(L_{\pm1})=0.$ \par

For any $n\in\Z$, one can write $D(L_n)$, $D(M_n)$, $D(N_n)$ and
$D(Y_{n-1/2})$ as follows:

\begin{eqnarray*}
&&D(L_n)=\mbox{$\sum\limits_{i\in\Z}$}(a_{1,i}^{(n)}L_i\!\otimes\!
L_{n-i}\!+\!a_{2,i}^{(n)}L_i\!\otimes\!M_{n-i}\!+\!a^{(n)}_{3,i}M_i\!\otimes\!
L_{n-i}\!+\!a_{4,i}^{(n)}L_i\!\otimes\!N_{n-i}\!+\!a_{5,i}^{(n)}N_i\!\otimes\!L_{n-i}\!\\
&&\ \ \ \ \ \ \
+a_{6,i}^{(n)}\!M_i\!\otimes\!M_{n-i}\!+\!a_{7,i}^{(n)}M_i\!\otimes\!N_{n-i}\!+\!a_{8,i}^{(n)}N_i\!\otimes\!M_{n-i}\!
+\!a_{9,i}^{(n)}N_i\!\otimes\!N_{n-i}\!+\!a_{10,i}^{(n)}
Y_{i-1/2}\!\otimes\!Y_{n-i+1/2}),\\
&&D(M_n)=\mbox{$\sum\limits_{i\in\Z}$}(b_{1,i}^{(n)}L_i\!\otimes\!
L_{n-i}\!+\!b_{2,i}^{(n)}L_i\!\otimes\!M_{n-i}\!+\!b^{(n)}_{3,i}M_i\!\otimes\!
L_{n-i}\!+\!b_{4,i}^{(n)}L_i\!\otimes\!N_{n-i}\!+\!b_{5,i}^{(n)}N_i\!\otimes\!L_{n-i}\!\\
&&\ \ \ \ \ \ \
+\!b_{6,i}^{(n)}\!M_i\!\otimes\!M_{n-i}\!+\!b_{7,i}^{(n)}M_i\!\otimes\!N_{n-i}\!+\!b_{8,i}^{(n)}N_i\!\otimes\!M_{n-i}\!
+\!b_{9,i}^{(n)}N_i\!\otimes\!N_{n-i}\!+\!b_{10,i}^{(n)}
Y_{i-1/2}\!\otimes\!Y_{n-i+1/2}),\\
&&D(N_n)=\mbox{$\sum\limits_{i\in\Z}$}(d_{1,i}^{(n)}L_i\!\otimes\!
L_{n-i}\!+\!d_{2,i}^{(n)}L_i\!\otimes\!M_{n-i}\!+\!d^{(n)}_{3,i}M_i\!\otimes\!
L_{n-i}\!+\!d_{4,i}^{(n)}L_i\!\otimes\!N_{n-i}\!+\!d_{5,i}^{(n)}N_i\!\otimes\!L_{n-i}\!\\
&&\ \ \ \ \ \ \
+\!d_{6,i}^{(n)}\!M_i\!\otimes\!M_{n-i}\!+\!d_{7,i}^{(n)}\!M_i\!\otimes\!N_{n-i}\!+\!d_{8,i}^{(n)}N_i\!\otimes\!M_{n-i}\!
+\!d_{9,i}^{(n)}N_i\!\otimes\!N_{n-i}\!+\!d_{10,i}^{(n)}
Y_{i-1/2}\!\otimes\!Y_{n-i+1/2}),\\
&&D(Y_{n-1/2})=\mbox{$\sum\limits_{i\in\Z}$}(\a_{n,i}L_i\!\otimes\!
Y_{n-1/2-i}\!+\!\a^\dag_{n,i}Y_{i-1/2}\!\otimes\!
L_{n-i}\!+\!\b_{n,i}M_i\!\otimes\!
Y_{n-1/2-i}\!+\!\b^\dag_{n,i}Y_{i-1/2}\!\otimes\!M_{n-i}\\&& \ \ \ \
\ \ \ \ \ \ \ \ \ \ +\gamma_{n,i}N_i\otimes Y_{n-1/2-i}
+\gamma_{n,i}^\dag Y_{i-1/2}\otimes N_{n-i}).
\end{eqnarray*}
Note that all the sums are finite. For any $n\in\Z$, one can easily
get the following identities by (\ref{LB}):
\begin{eqnarray*}
&&L_1\cdot(N_n\otimes N_{-n})=nN_{n+1}\otimes N_{-n}-nN_n\otimes
N_{1-n},\\
&&L_1\cdot(M_n\otimes N_{-n})=nM_{n+1}\otimes N_{-n}-nM_n\otimes
N_{1-n},\\
&&L_1\cdot(N_n\otimes M_{-n})=nN_{n+1}\otimes M_{-n}-nN_n\otimes
M_{1-n},\\
&&L_1\cdot(M_n\otimes M_{-n})=nM_{n+1}\otimes
M_{-n}-nM_n\otimes M_{1-n},\\
&&L_1\cdot(L_n\otimes N_{-n})=(n-1)L_{n+1}\otimes N_{-n}-nL_n\otimes
N_{1-n},\\
&&L_1\cdot(N_n\otimes L_{-n})=nN_{n+1}\otimes
L_{-n}-(1+n)N_n\otimes L_{1-n},\\
&&L_1\cdot(L_n\otimes M_{-n})=(n-1)L_{n+1}\otimes
M_{-n}-nL_n\otimes M_{1-n},\\
&&L_1\cdot(M_n\otimes L_{-n})=nM_{n+1}\otimes
L_{-n}-(1+n)M_n\otimes L_{1-n},\\
&&L_1\cdot(L_n\otimes L_{-n})=(n-1)L_{n+1}\otimes
L_{-n}-(1+n)L_n\otimes L_{1-n},\\
&&L_1\cdot(Y_{n-1/2}\otimes Y_{1/2-n})=(n-1)Y_{n+1/2}\otimes
 Y_{1/2-n}-nY_{n-1/2}\otimes Y_{3/2-n}.
\end{eqnarray*}
Let $Q_i=\rm{max}\{|p|\mid a^{(1)}_{i,p}\neq 0\}$ for $i=1,\! \cdots
,10.$ Applying $D-u_{\rm inn}$ to $L_1$, where $u$ is a proper
linear combination of $L_{p}\otimes L_{-p}$, $L_{p}\otimes M_{-p}$,
$M_{p}\otimes L_{-p}$, $M_{p}\otimes M_{-p}$, $L_{p}\otimes N_{-p}$,
$N_{p}\otimes L_{-p}$, $M_{p}\otimes N_{-p}$, $N_{p}\otimes M_{-p}$,
$N_{p}\otimes N_{-p}$ and $Y_{p-1/2}\otimes Y_{1/2-p}$ with
$p\in\Z$, and using induction on $\mbox{$\sum_{i=1}^{10}$} Q_i$, one
can safely suppose
\begin{eqnarray}\label{assumption1}
&&a_{1,i}^{(1)}=a_{6,j}^{(1)}=a_{7,j}^{(1)}=a_{8,j}^{(1)}=a_{9,j}^{(1)}=0,
\ \ \ \ \ {\rm for}\ i\neq-1,2,\,\ j\neq 0,1,\\ \label{assumption2}
&&a_{2,k}^{(1)}=a_{4,k}^{(1)}=a_{10,k}^{(1)}=a_{3,n}^{(1)}=a_{5,n}^{(1)}=0,
\ \ \ {\rm for}\  k\neq 0,2, \ \ n\neq\pm1.
\end{eqnarray}

Applying $D$ to $[\,L_{1},L_{-1}]=-2L_0$ and using $D(L_0)=0$, we
have \vspace{-5pt}
\begin{eqnarray*}\label{L1L-1}
\begin{aligned}
&\ \ \mbox{$\sum\limits_{p\in
\Z}$}\Big(\big((p-2)a_{1,p-1}^{(-1)}-(p+2)a_{1,p}^{(-1)}+(p-2)a_{1,p}^{(1)}-(p+2)a_{1,p+1}^{(1)}\big)L_p\otimes
L_{-p}\\[-10pt]
&\ \ \
+\big((p-2)a_{2,p-1}^{(-1)}-(p+1)a_{2,p}^{(-1)}+(p-1)a_{2,p}^{(1)}-(p+2)a_{2,p+1}^{(1)}\big)L_p\otimes
M_{-p}\\[-6pt]
&\ \ \
+\big((p-1)a_{3,p-1}^{(-1)}-(p+2)a_{3,p}^{(-1)}+(p-2)a_{3,p}^{(1)}-(p+1)a_{3,p+1}^{(1)}\big)M_p\otimes
L_{-p}\\[-6pt]
&\ \ \ +\big((p-2)a_{4,p-1}^{(-1)}-(p+1)a_{4,p}^{(-1)}+(p-1)a_{4,p}^{(1)}-(p+2)a_{4,p+1}^{(1)}\big)L_p\otimes N_{-p}\\[-6pt]
&\ \ \ +\big((p-1)a_{5,p-1}^{(-1)}-(p+2)a_{5,p}^{(-1)}+(p-2)a_{5,p}^{(1)}-(p+1)a_{5,p+1}^{(1)}\big)N_p\otimes L_{-p}\\[-6pt]
&\ \ \ +\big((p-1)a_{6,p-1}^{(-1)}-(p+1)a_{6,p}^{(-1)}+(p-1)a_{6,p}^{(1)}-(p+1)a_{6,p+1}^{(1)}\big)M_p\otimes M_{-p}\\[-6pt]
&\ \ \  +\big((p-1)a_{7,p-1}^{(-1)}-(p+1)a_{7,p}^{(-1)}+(p-1)a_{7,p}^{(1)}-(p+1)a_{7,p+1}^{(1)}\big)M_p\otimes N_{-p}\\[-6pt]
&\ \ \ +\big((p-1)a_{8,p-1}^{(-1)}-(p+1)a_{8,p}^{(-1)}+(p-1)a_{8,p}^{(1)}-(p+1)a_{8,p+1}^{(1)}\big)N_p\otimes M_{-p}\\[-6pt]
&\ \ \ +\big((p-1)a_{9,p-1}^{(-1)}-(p+1)a_{9,p}^{(-1)}+(p-1)a_{9,p}^{(1)}-(p+1)a_{9,p+1}^{(1)}\big)N_p\otimes N_{-p}\\[-6pt]
&\  \ \
+\big((p-2)a_{10,p-1}^{(-1)}-(p+1)a_{10,p}^{(-1)}+(p-2)a_{10,p}^{(1)}-(p+1)a_{10,p+1}^{(1)}\big)Y_{p-1/2}\otimes
Y_{1/2-p}\Big)=0.
\end{aligned}
\end{eqnarray*}
In particular, one has
\begin{eqnarray*}
(p-2)a_{1,p-1}^{(-1)}-(p+2)a_{1,p}^{(-1)}+(p-2)a_{1,p}^{(1)}-(p+2)a_{1,p+1}^{(1)}=0,\
\ \ \forall \ p\in\Z,
\end{eqnarray*}
which together with the fact that $\{p\in
\Z\,|\,a_{1,p}^{(-1)}\neq0\}$ is finite and (\ref{assumption1}),
forces
\begin{eqnarray}\label{bcd}
a_{1,p}^{(-1)}=3a_{1,-2}^{(-1)}+a^{(-1)}_{1,-1}+3a_{1,-1}^{(1)}
=a_{1,0}^{(-1)}+3a_{1,1}^{(-1)}+3a_{1,2}^{(1)}
=a_{1,-1}^{(-1)}+a_{1,0}^{(-1)}=0,
\end{eqnarray}
for  $p\in\Z\SM\{-2,0,\pm1\}$. Similarly, comparing the coefficients
of $L_p\otimes M_{-p}$, $M_p\otimes L_{-p}$ ,$L_p\otimes N_{-p}$,
$N_p\otimes L_{-p}$, $M_p\otimes M_{-p}$,  $M_p\otimes N_{-p}$,
$N_p\otimes M_{-p}$,  $N_p\otimes N_{-p}$  and $Y_{p-1/2}\otimes
Y_{1/2-p}$ and taking (\ref{assumption1}) and (\ref{assumption2})
into account, one has
\begin{eqnarray}
&&a_{2,0}^{(-1)}+2a_{2,1}^{(-1)}=a_{2,0}^{(-1)}+2a_{2,-1}^{(-1)}=a_{3,-1}^{(-1)}+2a_{3,0}^{(-1)}=a_{3,-1}^{(-1)}+2a_{3,-2}^{(-1)}=0,\label{bcd-3j}\\
&&a_{4,0}^{(-1)}+2a_{4,1}^{(-1)}=2a_{4,-1}^{(-1)}+a_{4,0}^{(-1)}+a_{4,0}^{(1)}=a_{5,-1}^{(-1)}+2a_{5,0}^{(-1)}=a_{5,-1}^{(-1)}+2a_{5,-2}^{(-1)}=0,\label{bcd-4j}\\
&&a_{i,-1}^{(-1)}+a_{i,0}^{(-1)}+a_{i,0}^{(1)}+a_{i,1}^{(1)}=2a_{10,-1}^{(-1)}\!+\!a_{10,0}^{(-1)}\!+\!2a_{10,0}^{(1)}=a_{10,0}^{(-1)}+2a_{10,1}^{(-1)}\!+\!2a_{10,2}^{(1)}=0,\label{bcd-2j}\\
&&a_{2,p}^{(1)}=a_{3,p}^{(1)}=a_{4,p}^{(1)}=a_{5,p}^{(1)}=a_{2,p_2}^{(-1)}=a_{3,p_3}^{(-1)}=a_{4,p_4}^{(-1)}=a_{5,p_5}^{(-1)}=a_{i,p_i}^{(-1)}=a_{10,p_{10}}^{(-1)}=0,\label{bcd-1j}
\end{eqnarray}
 for any $p\in\Z$, $p_2\in\Z\SM\{0,\pm 1\}$,
$p_3\in\Z\SM\{0,-1,-2\}$, $p_4\in\Z\SM\{0,\pm 1\}$,
$p_5\in\Z\SM\{0,-1,-2\}$, $p_i\in\Z\SM\{-1,0\}$ with  $i=6,7,8,9$
and $p_{10}\in\Z\SM\{0,\pm 1\}$. \par By (\ref{assumption1}) and
(\ref{assumption2}) as well as applying $D$ to $[L_1,N_0]=0$ and
$[L_1, M_0]=0$, respectively, we have
\begin{eqnarray}
&&a_{6,n}^{(1)}=a_{7,n}^{(1)}=a_{8,n}^{(1)}=a_{9,n}^{(1)}=a_{10,n}^{(1)}=0,\
\ \mbox{\ for\ \ all}\ n \in \Z,\label{f6}\\
&& d_{1,j_{1}}^{(0)}= d_{2,j_{2}}^{(0)}= d_{3,j_{3}}^{(0)}=
d_{4,j_{4}}^{(0)}=
 d_{5,j_{5}}^{(0)}= d_{i,j}^{(0)}= d_{10,p}^{(0)}=0,\label{f7}\\
 && b_{1,j_{1}}^{(0)}= b_{2,j_{2}}^{(0)}= b_{3,j_{3}}^{(0)}= b_{4,j_{4}}^{(0)}=
 b_{5,j_{5}}^{(0)}= b_{i,j}^{(0)}=b_{10,p}^{(0)}=0,\label{f8}\\
&&d_{1,0}^{(0)}+2d_{1,-1}^{(0)}=d_{1,0}^{(0)}+2d_{1,1}^{(0)}=d_{l,0}^{(0)}+d_{l,1}^{(0)}=d_{k,0}^{(0)}+d_{k,-1}^{(0)}=0,\label{f10}\\
&&b_{1,0}^{(0)}+2b_{1,-1}^{(0)}=b_{1,0}^{(0)}+2b_{1,1}^{(0)}=b_{l,0}^{(0)}+b_{l,1}^{(0)}=b_{k,0}^{(0)}+b_{k,-1}^{(0)}=0,\label{f99}
\end{eqnarray}
where $j_1\in\Z\SM\{0,\pm 1\}$, $j_2\in\Z\SM\{0,1\}$,
$j_3\in\Z\SM\{0,-1\}$, $j_4\in\Z\SM\{0,1\}$, $j_5\in\Z\SM\{0,-1\}$,
$j\in\Z\SM\{0\}$, $p\in\Z\SM\{0,1\}$, $i=6,\ 7,\ 8,\ 9$,  $l=2,\ 4,\
10$ and $k=3,\ 5$. Then it follows  from (\ref{assumption1}),
(\ref{bcd-1j}) and (\ref{f6}) that
\begin{eqnarray}\label{f3}
\ D(L_1)&=&\mbox{$\sum_{i\in\Z}$}\big(a_{1,i}^{(1)}L_i\otimes
L_{1-i} \big)=a_{1,-1}^{(1)}L_{-1}\otimes
L_2+a_{1,2}^{(1)}L_2\otimes L_{-1}.
\end{eqnarray}
Furthermore, applying $D$ to $[L_1, Y_{1/2}]=0$, we get $D(L_{1})=0$
from (\ref{f3}). Similarly,  applying $D$ to $[L_{-1},N_0]=0$ and
$[L_{-1}, M_0]=0$, respectively, we obtain from
(\ref{bcd-3j})-(\ref{bcd-1j}) and (\ref{f7})-(\ref{f99}) that
\begin{eqnarray}
&&a_{6,i}^{(-1)}=a_{7,i}^{(-1)}=a_{8,i}^{(-1)}=a_{9,i}^{(-1)}=a_{10,i}^{(-1)}=b_{4,i}^{(0)}=b_{5,i}^{(0)}=d_{4,i}^{(0)}=d_{5,i}^{(0)}=0,
\label{f11}\\
&&d_{2,1}^{(0)}-a_{2,0}^{(-1)}=d_{3,-1}^{(0)}-a_{3,-1}^{(-1)}=b_{2,1}^{(0)}+a_{4,0}^{(-1)}=b_{3,-1}^{(0)}+a_{5,0}^{(-1)}=0,\label{f12}
\end{eqnarray}
for all $i\in\Z$.  Set $u:=L_1\otimes M_{-1}-L_{0}\otimes M_0$.
Observe that $L_1\cdot u=0$. Substitute $D+a_{2,1}^{(-1)}u_{inn}$
into the expression of $D(L_{-1})$,  one can safely assume
$a_{2,1}^{(-1)}=0$, since such replacement would not affect the
expression of $D(L_1).$  Similarly, set $u^{(1)}:=M_{-1}\otimes
L_1-M_0\otimes L_0$, $u^{(2)}:=L_1\otimes N_{-1}-L_{0}\otimes N_0$,
and $u^{(3)}:=N_{-1}\otimes L_1-N_0\otimes L_0$, then replace $D$ by
$D+a_{3,0}^{(-1)}u^{(1)}_{inn}$, $D+a_{4,1}^{(-1)}u^{(2)}_{inn}$ and
$D+a_{5,0}^{(-1)}u^{(3)}_{inn}$ in turn, one can assume
$a_{3,0}^{(-1)}=a_{4,1}^{(-1)}=a_{5,0}^{(-1)}=0$. Hence we get
$$D(L_{-1})=a_{1,-2}^{(-1)}L_{-2}\otimes L_1+a_{1,-1}^{(-1)}L_{-1}\otimes L_0+a_{1,0}^{(-1)}L_{0}\otimes L_{-1}+a_{1,1}^{(-1)}L_{-1}\otimes L_2,$$
 by
(\ref{bcd})-(\ref{bcd-4j}). Finally, using $D([L_{-1},Y_{-1/2}])=0$
and (\ref{bcd}), one has $D(L_{-1})=0$.

 {\bf {Claim 3.}} $D(L_{\pm2})=D(N_1)=D(Y_{1/2})=0.$\par
It follows (\ref{f99}) and  (\ref{f11}) that
\begin{eqnarray}\begin{array}{l}\label{m}
D(M_0)=b_{1,-1}^{(0)}L_{-1}\otimes L_{1}+b_{1,0}^{(0)}L_{0}\otimes
L_{0} +b_{1,1}^{(0)}L_{1}\otimes L_{-1}+b_{6,0}^{(0)}M_0\otimes
M_0+b_{7,0}^{(0)}M_0\otimes N_0\\[8pt]\ \ \ \ \ \ \ \ \ \
+b_{8,0}^{(0)}N_0\otimes M_0+b_{9,0}^{(0)}N_0\otimes
N_0+b_{10,0}^{(0)}Y_{-1/2}\otimes
Y_{1/2}+b_{10,1}^{(0)}Y_{1/2}\otimes Y_{-1/2},
\end{array}
\end{eqnarray}
\vs{-0.7cm}
\begin{eqnarray}
\begin{array}{l}\label{n}
D(N_0)=d_{1,-1}^{(0)}L_{-1}\otimes L_{1}+d_{1,0}^{(0)}L_{0}\otimes
L_{0} +d_{1,1}^{(0)}L_{1}\otimes L_{-1}+d_{6,0}^{(0)}M_0\otimes
M_0+d_{7,0}^{(0)}M_0\otimes N_0\\[8pt]
\ \ \ \ \ \ \ \ \ \ \ +d_{8,0}^{(0)}N_0\otimes
M_0+b_{9,0}^{(0)}N_0\otimes N_0+d_{10,0}^{(0)}Y_{-1/2}\otimes
Y_{1/2}+d_{10,1}^{(0)}Y_{1/2}\otimes Y_{-1/2}.\ \ \ \ \ \ \ \
\end{array}
\end{eqnarray}
 Set $v^{(1)}:=M_0\otimes M_0$,  $v^{(2)}:=M_0\otimes N_0$,  $v^{(3)}:=N_0\otimes
 M_0$ and $v^{(4)}:=Y_{1/2}\otimes
Y_{-1/2}-Y_{-1/2}\otimes Y_{1/2}$. Observe that $L_{\pm 1}\cdot
v^{(i)}=0,$ but $N_0\cdot v^{(i)}\neq 0$ for $i=1,2,3,4$. Replacing
$D$ by $D-\frac14d_{6,0}^{(0)}v^{(1)}_{inn}$,
$D-\frac12d_{7,0}^{(0)}v^{(2)}_{inn}$,
$D-\frac12d_{8,0}^{(0)}v^{(3)}_{inn}$ and
$D-\frac12d_{10,0}^{(0)}v^{(4)}_{inn}$ in turn in (\ref{n}), one can
assume that
$d_{6,0}^{(0)}=d_{7,0}^{(0)}=d_{8,0}^{(0)}=d_{10,0}^{(0)}=0$. By
applying $D$ to $[N_0,N_1]=0$ and using (\ref{f10}), we have
$d_{1,-1}^{(0)}=d_{1,1}^{(0)}=d_{1,-1}^{(0)}=d_{10,1}^{(0)}=0$. Then
it follows from $D([N_0,M_0])=2D(M_0)$ and (\ref{m})-(\ref{n})that
$D(N_0)=0$ and
\begin{eqnarray}\label{f17}
D(M_0)=b_{7,0}^{(0)}M_0\otimes N_0+b_{8,0}^{(0)}N_0\otimes
M_0+b_{10,0}^{(0)}Y_{-1/2}\otimes
Y_{1/2}+b_{10,1}^{(0)}Y_{1/2}\otimes Y_{-1/2}\,.
\end{eqnarray}
Now considering $D([L_{\pm2},M_0])=0$ and $D([L_{\pm2},N_0])=0$, one
has
\begin{eqnarray}\label{f55}
D(L_2)=\mbox{$\sum_i $}a_{1,i}^{(2)}L_{i}\otimes L_{2-i}, \ \
\mbox{\ \ and}\ \ \ D(L_{-2})=\mbox{$\sum_i
$}a_{1,i}^{(-2)}L_{i}\otimes L_{-2-i}.
\end{eqnarray}
As a by-product, we also get $b_{10,0}^{(0)}=b_{10,1}^{(0)}=0$.
Replacing $D$ by $D+\frac12b_{8,0}^{(0)}(N_{0}\otimes N_0)_{inn}$ in
(\ref{f17}), one can assume $b_{8,0}^{(0)}=0$. Now we get from
(\ref{f17}) that
\begin{eqnarray}\label{f18}
D(M_0)=b_{7,0}^{(0)}M_0\otimes N_0.
\end{eqnarray}

Applying $D$ to $[L_1,L_{-2}]=-3L_{-1}$ and $[L_{-1},L_{2}]=3L_{1}$,
respectively, and using $D(L_{\pm 1})=0$, we have
\begin{eqnarray}
&&a_{1,p}^{(2)}=a_{1,q}^{(-2)}=2a_{1,1}^{(2)}+3a_{1,0}^{(2)}=a_{1,0}^{(2)}+4a_{1,-1}^{(2)}=a_{1,1}^{(2)}+a_{1,2}^{(2)}=a_{1,2}^{(2)}+4a_{1,3}^{(2)}=0,\label{f15}\\
&&
2a_{1,-1}^{(-2)}+3a_{1,0}^{(-2)}=a_{1,0}^{(-2)}+4a_{1,1}^{(-2)}=a_{1,-1}^{(-2)}+a_{1,-2}^{(-2)}=a_{1,-2}^{(-2)}+4a_{1,-3}^{(-2)}=0,
\label{f16}
\end{eqnarray}
where $p\in\Z\SM\{0,\pm1,2,3\}$ and  $q\in\Z\SM\{-3,-2,0,\pm1\}$.
Set $v:=L_{-1}\otimes L_1-2L_0\otimes L_0+L_1\otimes L_{-1}$ and
take $D-\frac14a_{1,0}^{(2)}v_{inn}$ in place of $D$ in the first
equation of (\ref{f55}), one can assume $a_{1,0}^{(2)}=0$. Then it
follows (\ref{f55}) and (\ref{f15}) that $D(L_{2})=0$. Consequently,
one can easily get $D(L_{-2})=0$ by applying $D$ to
$[L_{-2},L_2]=4L_0$ and using (\ref{f16}).

\par Applying $D$ to $[M_0,Y_{1/2}]=0$ and $[M_0,Y_{-1/2}]=0$,
respectively, one has
\begin{eqnarray}\label{f19}
\gamma_{1,i}^\dag=\gamma_{0,i}^\dag=\gamma_{0,j}=\gamma_{0,j}=0,\
\forall \ i\in\Z ,\ j\in\Z\SM\{0\}.
 \end{eqnarray}
 Similarly, when $D$ is applied
to $[N_0,Y_{1/2}]=Y_{1/2}$ and $[N_0,Y_{-1/2}]=Y_{-1/2}$,
respectively, it follows
$$\beta_{1,i}=\beta_{1,i}^\dag=\beta_{0,i}=\beta_{0,i}^\dag=0,\ \ \ \mbox{for\ \ all}\ \  i\in\Z.$$
Using $D([L_1,Y_{1/2}])=0$, $D([L_{-1},Y_{-1/2}])=0$ and $D(L_{\pm
1})=0$, we have
\begin{eqnarray}
&&\alpha_{1,i}=\alpha_{1,i}^\dag=\alpha_{0,i}^\dag=\alpha_{0,j}=0,\
\forall\ i\in\Z\SM\{0,1\},
j\in\Z\SM\{0,-1\},\\
&&\alpha_{1,0}+\alpha_{1,1}=\alpha_{0,0}+\alpha_{0,-1}=\alpha_{1,0}^\dag+\alpha_{1,1}^\dag=\alpha_{0,0}^\dag+\alpha_{0,1}^\dag=0.
\end{eqnarray}
Applying $D$ to $[L_1,Y_{-1/2}]=-Y_{1/2}$, one has
$\alpha_{0,0}=\alpha_{1,1},$ $\alpha_{0,0}^\dag=\alpha_{1,0}^\dag$
and $\gamma_{0,0}=\gamma_{1,0}.$ Let $a=\alpha_{0,0}$,
$b=\alpha_{0,0}^\dag$ and $c=\gamma_{0,0}$\,.
 Hence, we can rewrite $D(Y_{\pm 1/2})$ as
follows:
\begin{eqnarray}
&&D(Y_{1/2})=-aL_0\otimes Y_{1/2}+aL_1\otimes Y_{-1/2}+b
Y_{-1/2}\otimes L_1-b Y_{1/2}\otimes L_0+cN_0\otimes Y_{1/2},\label{f88}\ \ \ \ \ \ \ \ \\
&&D(Y_{-1/2})=aL_0\otimes Y_{- 1/2}-aL_{-1}\otimes Y_{1/2}+b
Y_{-1/2}\otimes L_0-b Y_{1/2}\otimes L_{-1}+cN_0\otimes Y_{-1/2}.\ \
\  \ \ \ \ \ \ \label{f90}
\end{eqnarray}
When $D$ is applied  to $[L_2,Y_{-1/2}]=-\frac32Y_{3/2}$ and
$[L_{-2},Y_{3/2}]=\frac52Y_{-1/2}$, respectively, one has $a=b=0$,
since $D(L_{\pm2})=0$. It follows $c=d_{7,0}^{(0)}=0$ by applying
$D$ to $[Y_{-1/2},Y_{1/2}]=M_0$, which proves
$D(Y_{\pm1/2})=D(M_0)=0$ by (\ref{f18}), (\ref{f88}) and
(\ref{f90}).\par

 Now it is left to calculate  $D(N_1)$. Firstly, using $D([N_0,
 N_1])=0$ and $D(N_0)=0$, we have
\begin{eqnarray}\label{f20}
 d_{2,i}^{(1)}=d_{3,i}^{(1)}=d_{6,i}^{(1)}=d_{7,i}^{(1)}=d_{8,i}^{(1)}=d_{10,i}^{(1)}=0,\
 \ \forall \ \ i\in\Z.
 \end{eqnarray}
Then applying $D$ to $[L_{-1}, N_1]=N_0$ and using $D(N_0)=0$, we
obtain
\begin{eqnarray}\begin{array}{l}
d_{1,i_1}^{(1)}=d_{4,i_2}^{(1)}=d_{5,i_3}^{(1)}=d_{9,i_4}^{(1)}=d_{1,0}^{(1)}+d_{1,1}^{(1)}=d_{1,1}^{(1)}+3d_{1,2}^{(1)}=d_{1,0}^{(1)}+3d_{1,-1}^{(1)}=0,\\[6pt]
d_{4,0}^{(1)}+2d_{4,1}^{(1)}=d_{4,0}^{(1)}+2d_{4,-1}^{(1)}=d_{5,1}^{(1)}+2d_{5,0}^{(1)}=d_{5,1}^{(1)}+2d_{5,2}^{(1)}=d_{9,0}^{(1)}+d_{9,1}^{(1)}=0,\label{f21}
\end{array}
\end{eqnarray}
where $i_1\in\Z\SM\{0,\pm1,2\}$, $i_2\in\Z\SM\{0,\pm1\}$,
$i_3\in\Z\SM\{0,1,2\}$, and $i_4\in\Z\SM\{0,1\}$. Finally, by
applying $D$ to $[N_1,Y_{-1/2}]=Y_{1/2}$ and (\ref{f21}) as well as
$D(Y_{\pm1/2})=0$, we get
$$d_{1,i}^{(1)}=d_{4,i}^{(1)}=d_{5,i}^{(1)}=d_{9,i}^{(1)}=0,\ \ \ \
 \forall \ \ i\in\Z,$$
which together with (\ref{f20}), yields $D(N_1)=0$.  Hence, the
claim is proved, so is the theorem, since $\LL$ is generated by
$L_{\pm1}$, $L_{\pm2}$, $N_1$ and $Y_{1/2}$.\QED\vskip5pt

The following lemma is very useful to the main theorem in the paper.
\begin{lemm}\rm \label{lemma6}
Suppose $v\in\VV$ such that $x\cdot v\in {\rm Im}(\,\rm Id-\tau)$
for all $x\in\LL.$ Then $v\in {\rm Im}(\,\rm Id-\tau)$.
\end{lemm}
\ni{\it Proof.}\ \ First note that $\LL\cdot {\rm Im}(\,\rm
Id-\tau)\subset {\rm Im}(\,\rm Id-\tau).$  We shall show that after
several steps in each of which $v$ is replaced by $v-u$ for some
$u\in {\rm Im}(\,\rm Id-\tau)$, the zero element is obtained, which
leads us to the result. Write $v=\sum_{n\in\frac{1}{2}\Z}v_n.$
Obviously,
\begin{eqnarray}\label{eqrx}
v\in {\rm Im}(\,{\rm Id}-\tau)\ \,\Longleftrightarrow \ \,v_n\in
{\rm Im}(\,{\rm Id}-\tau),\ \ \forall\,\,n\in\frac{1}{2}\Z.
\end{eqnarray}
Then $\sum_{n\in\frac{1}{2}\Z}nv_n=L_0\cdot v\in {\rm Im}(\,\rm
Id-\tau)$. By (\ref{eqrx}), $nv_n\in {\rm Im}(\,\rm Id-\tau).$ In
particular, $v_{n}\in {\rm Im}(\,\rm Id-\tau)$ if $n\ne0$. Thus when
replacing $v$ by $v-\sum_{n\in\frac{1}{2}\Z^*}v_n$, one can suppose
$v=v_0\in\VV_0$. Write
\begin{eqnarray*}
v\!\!\!&=&\!\!\!\mbox{$\sum\limits_{i\in\Z}$}(a_{i}L_i\otimes
L_{-i}+b_{i}L_i\otimes M_{-i}+c_{i}M_i\otimes L_{-i}+d_{i}M_i\otimes
M_{-i}+e_iY_{i-1/2}\otimes Y_{1/2-i}\\
\!\!\!&+&\!\!\!f_{i}L_i\otimes N_{-i}+g_{i}N_i\otimes
L_{-i}+h_{i}N_i\otimes N_{-i}+k_{i}M_i\otimes N_{-i}+r_{i}N_i\otimes
M_{-i}).
\end{eqnarray*}
Since all the elements of the forms $E_i\otimes F_{-i}-F_{-i}\otimes
E_i$ and $Y_{i-1/2}\otimes Y_{1/2-i}-Y_{1/2-i}\otimes Y_{i-1/2}$ are
contained in ${\rm Im}(\,\rm Id-\tau),$ where $\{E_i,
F_i\}\subset\{L_i,M_i,N_i\}$ for all $i\in\Z$. Replacing $v$ by
$v-u$, where $u$ is a combination of some of these elements, we can
assume
\begin{eqnarray}\label{wpqr2}
c_i=g_i=r_i=0,\ \forall\ \,i\in\Z;\ \ a_{i},\ d_{i},h_i\ne
0\,\Longrightarrow\ \,i>0\ \mbox{ or }\ i=0;\ \ e_{i}\ne
0\,\Longrightarrow\ \,i>0.
\end{eqnarray}
Then $v$ can be rewritten as
\begin{eqnarray}
v&=&\!\!\mbox{$\sum\limits_{i\in\Z_{+}}$}(a_{i}L_i\otimes
L_{-i}+d_iM_i\otimes M_{-i}+h_iN_i\otimes N_{-i})\nonumber\\
&+&\!\mbox{$\sum\limits_{i\in\Z}$}(b_{i}L_i\otimes
M_{-i}+f_iL_i\otimes N_{-i}+k_iM_i\otimes N_{-i})
+\mbox{$\sum\limits_{i\in\Z_{>0}}e_{i}$}Y_{i-1/2}\otimes
Y_{1/2-i}.\label{sm1}
\end{eqnarray}
Assume $a_{p}\ne 0$ for some $p>0$. Choose $q>0$ such that $q\ne p$.
Then $L_{p+q}\otimes L_{-p}$ appears in $L_{q}\cdot v,$ but
(\ref{wpqr2}) implies the term $L_{-p}\otimes L_{p+q}$ does not
appear in $L_q\cdot v$, which contradicts the fact that $L_q\cdot
v\in {\rm Im}(\,\rm Id-\tau)$. Hence we get $a_i=0,\
\forall\,\,i\in\Z^*$. Similarly, one can suppose $d_i=h_i=0,\
\forall\,\,i\in\Z^*$ and $e_i=0$, $\forall\,\,i\in\Z$. Then
(\ref{sm1}) becomes
\begin{eqnarray}\label{sm2}
v=\mbox{$\sum\limits_{i\in\Z}$}(b_{i}L_i\otimes M_{-i}+f_iL_i\otimes
N_{-i}+k_iM_i\otimes N_{-i}) +a_{0}L_0\otimes L_{0}+d_{0}M_0\otimes
M_{0}+h_0N_0\otimes N_0\,.
\end{eqnarray}
By ${\rm Im}(\,\rm Id-\tau)\subset{\rm Ker}(\,\rm Id+\tau)$ and our
hypothesis $\LL\cdot v\subset{\rm Im}(\,\rm Id-\tau)$, we have
\begin{eqnarray*}
0=(\,{\rm Id}+\tau)M_0\cdot v =-4h_{0}(M_0\otimes N_{0}+N_0\otimes
M_{0}) -2\mbox{$\sum\limits_{i\in\Z}$}f_i(L_i\otimes
M_{-i}+M_{-i}\otimes L_i)-4\mbox{$\sum\limits_{i\in\Z}$}k_i(
M_{i}\otimes M_{-i})\,.
\end{eqnarray*}
Comparing the coefficients, one gets $ h_0=0$ and $ f_{i}=k_i=0\!, \
\forall\,\,i\in\Z.$ Similarly, $(\,{\rm Id}+\tau)(N_0\cdot v)=0$
implies $d_0=0$ and $b_i=0$ for all $i\in \Z$; and $(\,{\rm
Id}+\tau)(L_1\cdot v)=0$ leads to $a_0=0$. Then the lemma follows
from (\ref{sm2}). \QED\vskip5pt

By now we have enough in hand to classify the Lie bialgebra
structures on the extended Schr\"odinger-Virasoro Lie algebra. The
following theorem is the central result of the paper.
\begin{theo}\label{theo2}\vskip-3pt Let $(\LL ,[\cdot,\cdot])$ be
the extended Schr\"odinger-Virasoro Lie algebra. Then each Lie
bialgebra structure on $\LL$ is triangular coboundary.
\end{theo}
\noindent{\it Proof.}\rm\ Let $(\LL ,[\cdot,\cdot],\D)$ be a Lie
bialgebra structure on $\LL$.  Thanks to Theorem \ref{theo1}, there
exists some $r\in\VV$ such that $\D=\D_r.$ By (\ref{cLie-s-s}),
$\mathrm{Im}(\D)\subset \mathrm{Im}(\,\rm Id-\tau)$. Then it follows
from Lemma \ref{lemma6} that $r\in\mathrm{Im}(\,\rm Id-\tau).$ But
(\ref{cLie-s-s}), Theorem \ref{theo}(iii) and  Lemma \ref{lemm3}
show that $\fc(r)=0$; as a result, $(\LL ,[\cdot,\cdot],\D)$ is a
triangular coboundary Lie bialgebra by Definition \ref{def2}.
 \QED

\end{document}